\documentclass{amsart}
\usepackage{amscd,amssymb,graphicx}
\author[Bodin]{Derek Bodin}
\address{
Derek Bodin\\
University of Minnesota\\
Minneapolis, MN 55455 }
\email{bodin@cs.umn.edu}
\author[DeCleene]{Chris DeCleene}
\address{
Chris DeCleene\\
University of Wisconsin-Eau Claire\\
Eau Claire, WI 54702-4004}
\email{cdecleene@gmail.com}
\author[Hager]{William Hager}
\address{
William Hager\\
University of Wisconsin\\
Eau Claire, WI 54702-4004}
\author[Otto]{Carolyn Otto}
\address{
Carolyn Otto\\
Rice University\\
Houston, TX 77005-1827} \email{cotto@rice.edu}
\author[Penkava]{Michael Penkava}
\address{
Michael Penkava\\
University of Wisconsin-Eau Claire\\
Eau Claire, WI 54702-4004} \email{penkavmr@uwec.edu}
\author[Phillipson]{Mitch Phillipson}
\address{
Mitch Phillipson\\
University of Wisconsin-Eau Claire\\
Eau Claire, WI 54702-4004}
\email{phillima@uwec.edu}
\author[Steinbach]{Ryan Steinbach}
\address{
Ryan Steinbach\\
University of Wisconsin-Madison\\
Madison, WI 53706-1796}
\email{rsteinbach@wisc.edu}
\author[Weber]{Eric Weber}
\address{
Eric Weber\\
University of Wisconsin-Eau Claire\\
Eau Claire, WI 54702-4004}
\email{webered@uwec.edu}
\thanks{Research of these authors was partially supported by grants
from the National Science Foundation and the University of
Wisconsin-Eau Claire.}
\input{xy}
\xyoption{all}

\newtheorem{thm}{Theorem}[section]

\theoremstyle{definition}

\def \ph{\varphi}

\def \ra{\rightarrow}

\def \hom{\mbox{\rm Hom}}
\def \ie{\hbox{\it i.e.}}

\def \tns{\otimes}

\def \mplus{+\cdots+}
\def \mcom{,\cdots,}
\def \k{\mbox{$\mathbb K$}}

\def \C{\mbox{$\mathbb C$}}
\def \Z{\mbox{$\mathbb Z$}}

\def\br#1#2{\lbrack#1,#2\rbrack}

\def\zt{\mbox{$\Z_2$}}

\def\inv{^{-1}}

\def\im{\operatorname{Im}}

\def\ainf{\mbox{$A_\infty$}}
\def\linf{\mbox{$L_\infty$}}
\def\and{\mbox{ \rm and }}
\def\T{\mathcal T}
\def\TV{\T(V)}
\def\TW{\mbox{$\T(W)$}}

\def\s#1{(-1)^{#1}}

\def\pha#1#2{\ph^{#1}_{#2}}

\def\psa#1#2{\psi^{#1}_{#2}}

\def\inv{^{-1}}

\newcommand{\mapsiso}{\tilde{\ra}}

\newcommand{\qbar}{\mbox{$\overline{\mathbb Q}$}}

\def \Cn#1#2{C^{#1}_{#2}}

\def\dec{\operatorname{Ch}}
\def\Ch{\operatorname{Ch}}
\def\Pn#1#2{P^{#1}_{#2}}
\def\Sn#1#2{S^{#1}_{#2}}
\def\la{\lambda}

\begin{document}
\setlength{\multlinegap}{0pt}
\title[The moduli space of complex two dimensional associative algebras]
{The moduli space of 2-dimensional associative algebras}%

\address{}%
\email{}%

\thanks{}%

\date{\today}
\begin{abstract}
In this paper, we study  moduli spaces of 2-dimensional 
complex associative algebras.
We give a complete calculation of the cohomology of every
element in the moduli space, as well as compute their versal deformations.
\end{abstract}
\maketitle

\nocite{ps1,ps2,mp1,pv1}

\section{Introduction}
A constructive approach to the computation of versal deformations of a large class of
algebraic algebras, including infinity algebras, was developed in
\cite{fp1}, and a procedure for computing examples was applied
in \cite{fp2} to construct some new examples of infinity
algebras, as well as to construct their miniversal deformations. Even though
infinity algebras have recently emerged as an important concept in
the area of mathematical physics called string theory, few finite dimensional
examples of such algebras have been studied, and their
versal deformations have not been constructed. The  examples in this paper
are the first step in the study of moduli spaces of low dimensional \ainf\
algebras.

In \cite{fp3,Ott-Pen}, the  moduli space of 3-dimensional Lie algebras was studied, and the miniversal
deformations of the Lie algebras played an important role in understanding how the moduli space was
glued together.  In \cite{fp8}, the moduli space of 4-dimensional Lie algebras was studied in the same manner.
The new perspective, using miniversal deformations to study the moduli spaces of Lie algebras, led to a more
complete picture of these moduli spaces.  The intention in this paper is to undertake a similar study for
2-dimensional associative algebras.

One method of constructing the moduli space of ordinary associative algebras in
dimension 2 is to consider extensions of a 1-dimensional associative algebra by
a 1-dimensional associative algebra. This is possible because there are no
simple 2-dimensional associative algebras, by a theorem of Wedderburn, so all
such algebras have an ideal, and therefore arise as extensions.  

In this paper, we will give a complete description of the moduli space of
2-dimensional complex associative
algebras, including a computation of a miniversal deformation
of each of these algebras.  From the miniversal deformations, a decomposition
of the moduli space into strata is obtained, with the only connections between
strata given by jump deformations. In the 2-dimensional case, the
description is simple, because each of the strata consists of a single point, so
the only interesting information is given by the jump deformations.

The versal deformation of an associative algebra depends only on the second and
third Hochschild cohomology groups.  However, we give a complete calculation of
the cohomology for each of the algebras. What makes the study of associative
algebras of low dimension much more complicated than the corresponding study of
low dimensional Lie algebras is that while for a Lie algebra, the $n$-th cohomology
group $H^n$ vanishes
for $n$ larger than the dimension of the vector space, in general, for an
associative algebra $H^n$ does not vanish.  Thus we had to develop arguments on
a case by case basis for each of the six distinct algebras.  In particular, one
of these algebras has an unusual pattern for the cohomology, which made its
computation rather nontrivial.

The main result of this paper is the complete description of the Hochschild
cohomology for
all $2$-dimensional associative algebras. It turns out that the calculation
of cohomology even for low dimensional associative algebras is a nontrivial
problem. To construct extensions of associative algebras to \ainf\ algebras, it
is necessary to have a complete description of the cohomology in all degrees,
not just $H^2$ and $H^3$, which are needed for the deformation theory of these
algebras as associative algebras. What we compute in this paper is the first
step in constructing $2$-dimensional \ainf\ algebras. These results may be of
interest on their own, especially as an indication of the difficulty which
occurs in computing the deformation theory of associative algebras, even in low
dimension.

\section{Preliminaries}
Suppose that $V$ is a 2-dimensional vector space, defined over a field $\k$
whose characteristic is not 2 or 3, equipped with an associative multiplication
structure $m:V\tns V\ra V$. The associativity relation can be given in the form
\begin{equation*}
m\circ(m\tns 1)=m\circ(1\tns m).
\end{equation*}

When the space $V$ is \zt-graded, there is no difference in the relation of
associativity, but only even maps $m$ are allowed, so the set of associative
algebra structures depends on the \zt-grading in this way.

The notion of \emph{equivalence} of associative algebra structures is given as
follows. If $g$ is a linear automorphism of $V$, then define
\begin{equation*}
g^*(m)=g\inv\circ m\circ (g\tns g).
\end{equation*}
Two algebra structures $m$ and $m'$ are equivalent if there is an automorphism
$g$ such that $m'=g^*(m)$. The set of equivalence classes of algebra structures
on $V$ is called the \emph{moduli space} of associative algebras on $V$.

When $V$ is \zt-graded, we require that $g$ be an even map. Thus the set of
equivalence classes of \zt-graded associative algebra structures will be
different than the set of equivalence classes of associative algebra structures
on the same space, ignoring the grading. Because the set of equivalences is
more restricted in the \zt-graded case, two algebra structures which are
equivalent as ungraded algebra structures may not be equivalent as \zt-graded
algebra structures. There is a map between the moduli space of \zt-graded
algebra structures on $V$ and the space of all algebra structures on $V$. In
general, this map will be neither injective nor surjective.

\emph{Hochschild cohomology} was introduced in \cite{hoch}, and used to
classify infinitesimal deformations of associative algebras. Suppose that
\begin{equation*}
m_t=m+t\ph,
\end{equation*}
is an infinitesimal deformation of $m$.  By this we mean that the structure
$m_t$ is associative up to first order. From an algebraic point of view, this
means that we assume that $t^2=0$, and then check whether associativity holds.
It is not difficult to show that is equivalent to the following.
\begin{equation*}
a\ph(b,c)-\ph(ab,c)+\ph(a,bc)-\ph(a,b)c=0,
\end{equation*}
where, for simplicity, we denote $m(a,b)=ab$. Moreover, if we let
\begin{equation*}
g_t=I+t\lambda
\end{equation*} be an infinitesimal automorphism of $V$, where
$\lambda\in\hom(V,V)$, then it is easily checked that
\begin{equation*}
g_t^*(m)(a,b)=ab+t(a\lambda(b)-\lambda(ab)+\lambda(a)b).
\end{equation*}
 This naturally leads
to a definition of the Hochschild coboundary operator $D$ on $\hom(\TV,V)$ by
\begin{align*}
D(\ph)(a_0\mcom a_n)=&a_0\ph(a_1\mcom a_n)+\s{n+1}\ph(a_0\mcom a_{n-1})a_n\\&+\sum_{i=0}^{n-1}\s{n+1}\ph(a_0\mcom a_{i-1},a_ia_{i+1},a_{i+2}\mcom a_n)
.
\end{align*}
If we set $C^n(V)=\hom(V^n,V)$, then $D:C^n(V)\ra C^{n+1}(V)$. One obtains the
following classification theorem for infinitesimal deformations.
\begin{thm} The equivalence classes of infinitesimal deformations $m_t$ of an associative algebra structure $m$ under the action of the
group of infinitesimal automorphisms on the set of infinitesimal deformations
are classified by the Hochschild cohomology group
\begin{equation*}
H^2(m)=\ker(D:C^2(V)\ra C^3(V))/\im(D:C^1(V)\ra C^2(V)).
\end{equation*}
\end{thm}
When $V$ is \zt-graded, the  only modifications that are necessary are that
$\ph$ and $\lambda$ are required to be even maps, so we obtain that the
classification is given by $H^2_e(V)$, the even part of the Hochschild
cohomology group.

We wish to transform this classical viewpoint into the more modern viewpoint of
associative algebras as being given by codifferentials on a certain coalgebra.
To do this, we first introduce the \emph{parity reversion} $\Pi V$ of a
\zt-graded vector space $V$. If $V=V_e\oplus V_o$ is the decomposition of $V$
into its even and odd parts, then $W=\Pi V$ is the \zt-graded vector space
given by $W_e=V_o$ and $W_o=V_e$. In other words, $W$ is just the space $V$
with the parity of elements reversed.

Denote the tensor (co)-algebra of $W$ by $\TW=\bigoplus_{k=0}^\infty W^k$,
where $W^k$ is the $k$-th tensor power of $W$ and $W^0=\k$. For brevity, the
element in $W^k$ given by the tensor product of the elements $w_i$ in $W$ will
be denoted by $w_1\cdots w_k$. The coalgebra structure on $\TW$ is  given by
\begin{equation*}
\Delta(w_1\cdots w_n)=\sum_{i=0}^n w_1\cdots w_i\tns w_{i+1}\cdots w_n.
\end{equation*}
Define $d:W^2\ra W$ by $d=\pi\circ m\circ (\pi\inv\tns\pi\inv)$, where
$\pi:V\ra W$ is the identity map, which is odd, because it reverses the parity
of elements. Note that $d$ is an odd map. The space $C(W)=\hom(\TW,W)$ is
naturally identifiable with the space of coderivations of $\TW$.  In fact, if
$\ph\in C^k(W)=\hom(W^k,W)$, then $\ph$ is extended to a coderivation of $\TW$
by
\begin{equation*}
\ph(w_1\cdots w_n)=
\sum_{i=0}^{n-k}\s{(w_1\mplus w_i)\ph}w_1\cdots
 w_i\ph(w_{i+1}\cdots w_{i+k})w_{i+k+1}\cdots w_n.
\end{equation*}

The space of coderivations of $\TW$ is equipped with a \zt-graded Lie algebra
structure given by
\begin{equation*}
[\ph,\psi]=\ph\circ\psi-\s{\ph\psi}\psi\circ\ph.
\end{equation*}
The reason that it is more convenient to work with the structure $d$ on $W$
rather than $m$ on $V$ is that the condition of associativity for $m$
translates into the codifferential property $[d,d]=0$.  Moreover, the
Hochschild coboundary operation translates into the coboundary operator $D$ on
$C(W)$, given by
\begin{equation*}
D(\ph)=[d,\ph].
\end{equation*}
This point of view on Hochschild cohomology first appeared in \cite{sta4}.  The
fact that the space of Hochschild cochains is equipped with a graded Lie
algebra structure was noticed much earlier \cite{gers,gers1,gers2,gers3,gers4}.

For notational purposes, we introduce a basis of $C^n(W)$ as follows.  Suppose
that $W=\langle w_1\mcom w_m\rangle$. Then if $I=(i_1\mcom i_n)$ is a
\emph{multi-index}, where $1\le i_k\le m$, denote $w_I=w_{i_1}\cdots w_{i_n}$.
Define $\ph^{I}_i\in C^n(W)$ by
\begin{equation*}
\ph^I_i(w_J)=\delta^I_Jw_i,
\end{equation*}
where $\delta^I_J$ is the Kronecker delta symbol. In order to emphasize the
parity of the element, we will denote $\ph^I_i$ by $\psi^I_i$ when it is an odd
coderivation.

For a multi-index $I=(i_1\mcom i_k)$, denote its \emph{length}  by $\ell(I)=k$.  If
$K$ and $L$ are multi-indices, then denote $KL=(k_1\mcom k_{|K|},l_l\mcom
l_{|L|})$.  Then
\begin{align*}
(\ph^I_i\circ\ph^J_j)(w_K)&=
\sum_{K_1K_2K_3=K}\s{w_{K_1}\ph^J_j} \ph^I_i(w_{K_1},\ph^J_j(w_{K_2}), w_{K_3})
\\&=
\sum_{K_1K_2K_3=K}\s{w_{K_1}\ph^J_j}\delta^I_{K_1jK_3}\delta^J_{K_2}w_i,
\end{align*}
from which it follows that
\begin{equation}\label{braform}
\ph^I_i\circ\ph^J_j=\sum_{k=1}^{\ell(I)}\s{(w_{i_1}\mplus w_{i_{k-1}})\ph^J_j}
\delta^k_j
\ph^{(I,J,k)}_i,
\end{equation}
where $(I,J,k)$ is given by inserting $J$ into $I$ in place of the $k$-th
element of $I$; \ie, $(I,J,k)=(i_1\mcom i_{k-1},j_1\mcom j_{|J|},i_{k+1}\mcom
i_{\ell(I)})$.

Let us recast the notion of an infinitesimal deformation in terms of the
language of coderivations.  We say that
\begin{equation*}
d_t=d+t\psi
\end{equation*}
is a deformation of the codifferential $d$ precisely when $[d_t,d_t]=0 \mod t^2$.
This condition immediately reduces to the cocycle condition $D(\psi)=0$.  Note
that we require $d_t$ to be odd, so that $\psi$ must be an odd coderivation.
One can introduce a more general idea of parameters, allowing both even and odd
parameters, in which case even coderivations play an equal role, but we will
not adopt that point of view in this paper.

For associative algebras, we require that $d$ and $\psi$ lie in $\hom(W^2,W)$.
This notion naturally generalizes to considering $d$ simply to be an arbitrary
odd codifferential, in which case we would obtain an \ainf\ algebra, a natural
generalization of an associative algebra.

\section{Associative algebra structures on an $0|2$ vector space}
Suppose that $W=\langle 1,2 \rangle$,  where both $1,2$ are odd elements. Then $C^n=\langle \ph^I_i,\ell(I)=n\rangle$ has dimension $\dim C^n=2^{n+1}$. For later convenience, we decompose $C^n$ as follows. Let
\begin{align*}
\Cn{n}{1}=\langle \ph^I_1, \ell(I)=n\rangle\\
\Cn{n}{2}=\langle \ph^I_2,\ell(I)=n\rangle.
\end{align*}
Then $C^n=\Cn{n}{1}\oplus \Cn{n}{2}$. Moreover $\dim \Cn{n}{1}=\dim \Cn{n}{2}=2^n$.

The matrix of a generic odd coderivation is
$A=\left[ \begin {array}{cccc} a_{{1,1}}&a_{{1,2}}&a_{{1,3}}&a_{{1,4}}\\\noalign{\medskip}a_{{2,1}}&a_{{2,2}}&a_{{2,3}}&a_{{2,4}}
\end {array} \right]
$
The codifferential condition $[d,d]=0$ gives 8 solutions,
\begin{enumerate}
\item $a_{{1,2}}=a_{2,3}=a_{2,1}=a_{1,4}=a_{1,1}=a_{2,2}=0,\quad a_{{1,3}}=a_{{2,4}}$
\item $a_{{1,3}}=a_{{2,3}}=a_{{2,1}}=a_{{1,4}}=a_{{1,1}}=a_{{2,2}}=0,\quad a_{{1,2}}=a_{{2,4}}$
\item $a_{{1,2}}=a_{{1,3}}=a_{{2,3}}=a_{{2,1}}=a_{{1,4}}=a_{{2,2}}=0,\quad a_{{1,1}}=a_{{2,4}}$
\item $a_{{2,3}}=a_{{1,4}}=a_{{2,2}}=0,\quad a_{{1,3}}=a_{2,4}=a_{{1,2}}$
\item $a_{{1,3}}=a_{{2,1}}=a_{{1,4}}=a_{{2,2}}=0,\quad a_{{2,4}}=a_{{1,2}},\quad a_{{1,1}}=a_{{2,3}}$
\item $a_{{1,2}}=a_{{2,3}}=a_{{2,1}}=a_{{1,4}}=0,\quad a_{{2,4}}=a_{{1,3}},\quad a_{{1,1}}=a_{{2,2}}$
\item $a_{{1,2}}=a_{{1,3}}=a_{{1,4}}=0,\quad a_{{2,2}}=a_{{2,3}},\quad a_{{1,1}}=-{\frac {a_{{2,1}}a_{{2,4}}-{a_{{2,3}}}^{2}}{a_{{2,3}}}}$
\item $a_{{2,2}}=a_{{2,3}},\quad a_{{1,2}}=a_{{1,3}},\quad a_{{2,1}}={\frac {a_{{1,3}}a_{{2,3}}}{a_{{1,4}}}},
\quad a_{{1,1}}={\frac {{a_{{1,3}}}^{2}-a_{{1,3}}a_{{2,4}}+a_{{1,4}}a_{{2,3}}}{a_{{1,4}}}}$.
\end{enumerate}

Each of these solutions corresponds to at least one nonequivalent codifferential listed below in the following manner:
\begin{itemize}
\item Number (1) corresponds to $d_4$. This correspondence can be seen by setting $a_{2,4}=1$ and applying an automorphism that interchanges $w_1$ with $w_2$.
\item Number (2) corresponds to $d_3$, for reasons identical to those above.
\item Number (3) corresponds to $d_1$. This correspondence can be seen by setting $a_{2,4}=1$.
\item Number (4) corresponds to $d_1$,$d_2$, and $d^6$. For $d_2$ we set $a_{2,4}=0$. For $d_6$ we set $a_{2,2}=0$ and $a_{2,4}=1$ and then apply an automorphism similar to that mentioned with regards to number (1). Though the arguments are similar for $d_1$, the automorphisms are slightly more complicated.
\item Number (5) also corresponds to $d_3$. This correspondence can be seen by setting $a_{1,1}=1$ and $a_{2,4}=0$.
\item Number (6) corresponds to $d_4.$ This correspondence is identical to that above.
\item Number (7) also corresponds to $d_1$,$d_2$, and $d_6$. For $d_6$, we set $a_{2,3}=1$ and everything else to 0. The other two correspondences have slightly more complicated automorphisms.
\item Number (8) has the least restrictions and so corresponds to $d_1$,$d_2$,$d_5$, and $d_6$.
\end{itemize}

Nonequivalent Codifferentials:
\begin{align*}
d_1&=\psi_2^{22}+\psi_1^{11}\\
d_2&=\psi_2^{22}\\
d_3&=\psi_2^{22}+\psi^{12}_1\\
d_4&=\psi_2^{22}+\psi_1^{21}\\
d_5&=\psi_2^{22}+\psi_1^{12}+\psi_1^{21}\\
d_6&=\psi_1^{22}.
\end{align*}
 Note that if we define $D(\ph)=\br {d^*}\ph$, where $d^*$ is one of the above codifferentials, then $D^2=0$, so the \emph{coboundary operator} $D$ determines a differential on $C(W)$. Since $d^*\in C^2$, $D(C^k)\subseteq C^{k+1}$, and we can define the $k$-th cohomology $H^k(d^*)$ of $d^*$ by
\begin{equation*}
H^k(d^*)=\ker(d^*:C^k\ra C^{k+1})/\im(d^*:C^{k-1}\ra C^k).
\end{equation*}

The cohomology of these codifferentials is given in Table \ref{coho table}
below.

\begin{table}[h!]
\begin{center}
\begin{tabular}{lccccc}
Codifferential&$H^0$&$H^2$&$H^1$&$H^3$&$H^4$\\ \hline \\
$d_1=\psi_1^{11}+\psi_2^{22}$&2&0&0&0&0\\
$d_2=\psi_2^{22}$&2&1&1&1&1\\
$d_3=\psi_2^{22}+\psi^{12}_1$&0&0&0&0&0\\
$d_4=\psi_2^{22}+\psi_1^{21}$&0&0&0&0&0\\
$d_5=\psi_2^{22}+\psi_1^{12}+\psi_1^{21}$&2&1&1&1&1\\
$d_6=\psi_1^{22}$&2&2&2&2&2\\\\ \hline
\end{tabular}
\end{center}
\label{coho table}
\caption{Cohomology of the six codifferentials on a $0|1$-dimensional space}
\end{table}

\section{Elements of the moduli space}
In this section we give a complete description of both the cohomology and the
multiplication structure generated by each codifferential. For a complete proof
of the cohomological structure see the next section. Let us suppose that
$V=\langle x,\theta\rangle$, where $x,\theta$ are even, and that $W=\Pi V=\langle 1,2\rangle$,
where $\pi(x)=2$ and $\pi(\theta)=1$. Let $m=\pi\circ d\circ (\pi\inv\tns\pi\inv)$.
Then $m$ is an associative algebra structure on $V$, corresponding to the codifferential $d$.
For each of the codifferentials, we give the multiplication structure $m$ on $V$.
\begin{equation*}
\begin{array}{lllll}
d_1&x^2=x&x\theta=0&\theta x=0&\theta^2=\theta\\
d_2&x^2=0&x\theta=0&\theta x=0&\theta^2=\theta\\
d_3&x^2=0&x\theta=x&\theta x=0&\theta^2=\theta\\
d_4&x^2=0&x\theta=0&\theta x=x&\theta^2=\theta\\
d_5&x^2=0&x\theta=x&\theta x=x&\theta^2=\theta\\
d_6&x^2=0&x\theta=0&\theta x=0&\theta^2=x\\
\end{array}
\end{equation*}
Of these algebras, $d_1$, $d_2$, $d_5$ and $d_6$ are commutative;
and $d_1$ and $d_5$ are unital, with unit $\theta$. In the algebras $d_1$ and
$d_5,$ $\theta$ generates a nontrivial proper ideal, while $x$ generates a nontrivial
proper ideal in every algebra. The algebra $d_1$ is the unique semisimple 2-dimensional
algebra, which is the direct sum of two copies of the 1-dimensional simple algebra $\C$.

The algebras $d_2$, $d_3$, $d_4$ and $d_5$ are all extensions of the simple $1$-dimensional
associative algebra (whose structure is just the associative algebra structure of $\C$). In
fact, they fit a certain pattern of extensions. The algebras $d_3$ and $d_4$ are opposite
algebras, and they are rigid in the cohomological sense. These two rigid algebras are just
the first in a sequence of rigid extensions of the $1$-dimensional simple algebra.

The algebra $d_5$ is the unique
extension of the simple $1$-dimensional algebra by the trivial $1$-dimensional
algebra as a unital algebra. The algebra $d_2$ is just the direct sum of the
trivial $1$-dimensional algebra and the simple $1$-dimensional algebra.

Finally, the algebra $d_6$ is an extension of the trivial $1$-dimensional algebra by
the trivial $1$-dimensional algebra, and as a consequence, it is a nilpotent algebra.
By nilpotent algebra, we mean an algebra such that a power of the algebra vanishes,
which in the finite dimensional case is equivalent to the fact that every element in
this algebra is nilpotent.

We did not use the method of extensions in calculating the nonequivalent codifferentials.
In this simple case, it is easy to solve the codifferential property $[d,d]=0$, which
gives a system of quadratic coefficients, and study the action of the group of linear
automorphisms of the underlying vector space, to arrive at the six codifferentials.
However, calculating this space by extensions reveals more of its properties, and also
gives a natural manner of organizing the codifferentials.

\section{Calculating the cohomology}
The cohomology of the codifferentials is given in Table \ref{coho table} above.
With the exception of $d_6$, the pattern of cohomology is easily deduced from
the information in the table.

For later use, we define  the following operator on $C(W)$. If $I$ is a
multi-index with $i_k\in\{1,2\}$, with $\ell(I)=m$, then define $\lambda^I:C^k\ra
C^{k+m}$ by $\lambda^I\pha Jj=\pha{IJ}j$. Note that the parity of $\lambda^I$
is the same as the parity of $I$. We abbreviate $\lambda^{\{1\}}=\lambda^1$.

We give a computation of the cohomology of the codifferentials on a case by
case basis.

\begin{thm}
 Suppose that a coboundary operator $D:C^n\ra C^{n+1}$ decomposes as $D=D'+D''$, given by the following diagram
 $$\xymatrix@C=13pt@R=15pt{
 \Cn{n}{a} \ar[d]^{D'} \ar[rd]^{D''} \\
 \Cn{n+1}{a} \ar[d]^{D'} \ar[rd]^{D''} \ar @{} [r]|{\oplus} &
 \Cn{n+1}{a+1} \ar[d]^{D'} \ar[rd]^{D''} \\
 \Cn{n+2}{a} \ar[d]^{D'} \ar[rd]^{D''} \ar @{} [r]|{\oplus} &
 \Cn{n+2}{a+1} \ar[d]^{D'} \ar[rd]^{D''}\ar@{}[r]|{\oplus} &
 \Cn{n+2}{a+2} \ar[d]^{D'} \ar[rd]^{D''} \\
 &&&&
  }
  $$
  for $a\le k\le n$ where
  $\Cn{n}{} =\Cn{n}1\oplus\dots\oplus\Cn{n}a\oplus\dots\oplus\Cn{n}{n}$,
  such that $D''$ is injective when $k=a$, and  $H(D'')=0$. Then $H(D)=0$ on the subcomplex
$\Cn{n}{k}$ for $k\ge a$.
\label{thm-arg}
\end{thm}
This result is well known. A proof is given in \cite{bdhoppsw1}.


\subsection{$d_1=\psi_2^{22}+\psi_1^{11}$}\label{01case3}
We begin by computing the coboundary operator with representatives from the $\Cn{n}{2}$
and $\Cn{n}{1}$ spaces,
\begin{align*}
D(\ph^I_1)=&\ph^{I1}_1+(-1)^{I+1}\ph^{1I}_1+(-1)^I\ph^I_1\psi^{22}_2+(-1)^I\ph^I_1\psi^{11}_1\\
D(\ph^I_2)=&\ph^{I2}_2+(-1)^{I+1}\ph^{2I}_2+(-1)^I\ph^I_2\psi^{22}_2+(-1)^I\ph^I_2\psi^{11}_1
\end{align*}
We decompose  the  $\Cn{n}{1}$ and $\Cn{n}{2}$ spaces as follows.
\begin{align*}
\Cn{n}{1,k}=&\langle\ph^{1^k2I}_1|\ell(I)=n-k-1\rangle & \Pn{n}1=&\langle\ph^{1^n}_1\rangle\\
\Cn{n}{2,k}=&\langle\ph^{2^k1I}_2|\ell(I)=n-k-1\rangle & \Pn{n}2=&\langle\ph^{2^n}_2\rangle
\end{align*}
We also decompose the coboundary operator $D$ as follows.
\begin{align*}
D=D'_1+D''_1:\Cn{n}{1,k}\ra&\Cn{n+1}{1,k}\oplus\Cn{n+1}{1,k+1} \\
D=D'_2+D''_2:\Cn{n}{2,k}\ra&\Cn{n+1}{2,k}\oplus\Cn{n+1}{2,k+1} \\
D:\Pn{n}a\ra&\Pn{n+1}a.
\end{align*}
By computation we see,
\begin{align*}
 D''_2(\ph^{2^k1I}_2)=&\left\{\begin{array}{cc}
                               \ph^{2^{k+1}1I}_2, & k\text{ is even}; \\
                               0, & k\text{ is odd}
                              \end{array}\right.
&D(\ph^{2^n}_2)=&\left\{\begin{array}{cc}
                               0, & n\text{ is even}; \\
                               \ph^{2^{n+1}}_2, & n\text{ is odd}
                              \end{array}\right. \\
 D''_1(\ph^{1^k2I}_1)=&\left\{\begin{array}{cc}
                               \ph^{1^{k+1}2I}_1, & k\text{ is even}; \\
                               0, & k\text{ is odd}
                              \end{array}\right.
& D(\ph^{1^n}_1)=&\left\{\begin{array}{cc}
                               0, & n\text{ is even}; \\
                               \ph^{1^{n+1}}_1, & n\text{ is odd}
                              \end{array}\right.
\end{align*}
Using Theorem \ref{thm-arg} and a direct computation of $H^0(d_1)$, we obtain that
$$H^n(d_1)=\left\{
        \begin{array}{ll}
          \langle\psi_2,\psi_1\rangle, & n=0; \\
          0, & n\geq 1.\\
        \end{array}
      \right.$$


\subsection{$d_2=\psi_2^{22}$}\label{01case1}
We will begin by computing the coboundary operator with representatives from the $\Cn{n}{2}$ and $\Cn{n}{1}$ spaces:
\begin{align*}
D(\ph^I_2) = &\ph^{I2}_2 + (-1)^{I+1}\ph^{2I}_2+(-1)^I\ph^I_2\psi^{22}_2 \\
D(\ph^I_1)=&(-1)^I\ph^I_1 \psi^{22}_2.
\end{align*}
We decompose our spaces as follows,
\begin{align*}
 \Cn{n}{a,k}=&\langle\ph^{2^k1I}_a\mid\ell(I)=n-k-1\rangle & \Pn{n}{a}=&\langle\ph^{2^n}_a\rangle.
\end{align*}
Using this decomposition we have the following maps,
\begin{align*}
 D=D'_a+D''_a:\Cn{n}{a,k}\ra&\Cn{n+1}{a,k}\oplus\Cn{n+1}{a,k+1} \\
 D:\Pn{n}a\ra&\Pn{n+1}a. \\
\end{align*}
By computation we see,
\begin{align*}
 D''_2(\ph^{2^k1I}_2)=&\left\{\begin{array}{cc}
                               \ph^{2^{k+1}1I}_2, & k\text{ is even}; \\
                               0, & k\text{ is odd}
                              \end{array}\right.
&D(\ph^{2^n}_2)=&\left\{\begin{array}{cc}
                               0, & n\text{ is even}; \\
                               \ph^{2^{n+1}}_2, & n\text{ is odd}
                              \end{array}\right. \\
 D''_1(\ph^{2^k1I}_1)=&\left\{\begin{array}{cc}
                               0, & k\text{ is even}; \\
                               \ph^{2^{k+1}1I}_1, & k\text{ is odd}
                              \end{array}\right.
& D(\ph^{2^n}_1)=&\left\{\begin{array}{cc}
                               0, & n\text{ is even}; \\
                               \ph^{2^{n+1}}_1, & n\text{ is odd}
                              \end{array}\right.
\end{align*}
Thus using Theorem \ref{thm-arg} we see that $H(D)=0$ on the subcomplexes $\Cn{n}{2,k}$ for $k\ge0$, $\Cn{n}{1,k}$ for $k\ge1$ and $\Pn{n}{a}$ for $n\ge1$. If we let $D_1=D'_1+D''_1$, then we see that $\lambda^1D_1=D_1\lambda^1$ which implies the following diagram,
$$\xymatrix@R=20mm@C=8mm{
H^n(C_1)\ar@1{-}[d] \ar@1{-}[r] & H^n(C_{1,0}) \ar@/^/[ld]^{\lambda^1} \\
H^{n+1}(C_{1,0}) &
}$$
This is true for $n\ge0$ if we define $\Cn0{1,0}=\Cn{0}{1}=\Pn01=\langle\ph_1\rangle$. Therefore we conclude,
$$H^n(d_2)=\left\{
        \begin{array}{ll}
          \langle\psi_2,\psi_1\rangle, & n=0; \\
          \langle\phi^{1^n}_1\rangle, & n\geq 1.\\
        \end{array}
      \right.$$


\subsection{$d_3=\psi_2^{22}+\psi_1^{12}$}\label{case2}
We begin by computing the coboundary operator with representatives from the $\Cn{n}{2}$ and $\Cn{n}{1}$ spaces,
\begin{align*}
D(\ph^I_2) =& \ph^{I2}_2 + (-1)^{I+1}\ph^{2I}_2+(-1)^I\ph^I_2\psi^{22}_2+(-1)^I\ph^I_2\psi^{12}_1+(-1)^{I+1}\ph^{1I}_1 \\
D(\ph^I_1) =& \ph^{I2}_1 + (-1)^I\ph^I_1\psi^{22}_2 + (-1)^I\ph^I_1\psi^{12}_1.
\end{align*}
Now we decompose $\Cn{n}1$ as follows,
\begin{align*}
 \Cn{n}{1,k}=&\langle\ph^{2^k1I}_1\mid\ell(I)=n-k-1\rangle,\qquad \Pn{n}{1}=\langle\ph^{2^n}_1\rangle \\
\end{align*}
Using this decomposition we have,
\begin{align*}
 D=D_2+D_1:\Cn{n}{2}\ra&\Cn{n+1}{2}\oplus\Cn{n+1}{1,0} \\
 D=D'+D'':\Cn{n}{1,k}\ra&\Cn{n+1,1}{k}\oplus\Cn{n+1}{1,k+1} \\
 D:\Pn{n}1\ra& \Pn{n+1}1
\end{align*}
Furthermore we see that $D:\Cn{n}{1,0}\ra\Cn{n+1}{1,0}$. Now we check how $D''$ acts
\begin{align*}
 D''(\ph^{2^k1I}_1)=\left\{\begin{array}{cc}
                            0 & k\text{ is even}; \\
                            \ph^{2^{k+1}1I}_1 & k\text{ is odd}
                           \end{array}\right.
\end{align*}
Thus by Theorem \ref{thm-arg} we see that the cohomology vanishes on the subcomplex $\Cn{n}{1,k},$ for $k\ge1$. Notice that $D_1:\Cn{n}{2}\mapsiso\Cn{n+1}{1,0}$ and that
\begin{align*}
 DD_1(\ph^I_2)=&\ph^{1I2}_1+(-1)^{I+1}\ph^{1I}_1\psi^{22}_2+(-1)^{I+1}\ph^{12I}_1+(-1)^I\lambda^1\ph^I_1\psi^{12}_1 \\
              =&(-1)^{n+1}D_1D_2(\ph^I_2)
\end{align*}
Let $a\ph\in\Cn{n}{2}$ for $a\in\C$ and $\xi\in\Cn{n}{1,0}$ be such that $a\ph+\xi\in\ker D$, then
\begin{align*}
 D_1(a\ph)+D(\xi)=&0 & D_2(a\ph)=&0.
\end{align*}
However the second equation follows from the first,
\begin{align*}
 0=&DD_1(a\ph)+D^1(\xi) \\
  =&(-1)^{n+1}D_1D_2(a\ph)
\end{align*}
and since $D_1$ is injective then $D_2(a\ph)=0$. Since $D_1$ is an isomorphism then $\xi=D_1(\eta)$ for $\eta\in\Cn{n-1}{2}$, thus
\begin{align*}
 0=&D_1(a\ph)+D(\xi) \\
  =&D_1(a\ph)+DD_1(\eta) \\
  =&D_1(a\ph+(-1)^{n+1}D_2(\eta)).
\end{align*}
This implies that $a\ph=(-1)^nD_2(\eta)$. Therefore we have shown $(-1)^na\ph+\xi=D(\eta)$ since $a$ is arbitrary we can absorb the $(-1)^n$, and we have
\begin{equation*}
 a\ph+\xi=D(\eta)
\end{equation*}

Next we consider is $\Pn{n}{1}$, but
\begin{align*}
 D(\ph^{2^n}_1)=&\left\{\begin{array}{ll}
           0, &n\text{ odd} \\
           \ph^{2^{n+1}}_1, &n\text{ even}
                  \end{array}
                  \right.
\end{align*}

We've now shown for all spaces except $\Cn{0}{2}$, but for this we calculate
\begin{align*}
 D(\ph_2)=-\ph^1_1\ne0
\end{align*}
So we discover
$$H^n(d_3)=0\quad\text{for all }n$$


\subsection{$d_4=\psi_2^{22}+\psi_1^{21}$}\label{01case4}
This is analogous to section \ref{case2}. We define the spaces as follows,
\begin{align*}
 \Cn{n}{2,k}=&\langle\ph^{I12^k}_2\mid\ell(I)=n-k-1\rangle,\quad &\Pn{n}2=&\langle\ph^{2^n}_2\rangle \\
 \Cn{n}{1,k}=&\langle\ph^{I12^k}_1\mid\ell(I)=n-k-1\rangle,\quad &\Pn{n}1=&\langle\ph^{2^n}_1\rangle
\end{align*}
The cohomology is given by,
$$H^n(d_4)=0\quad\text{for all }n$$


\subsection{$d_5=\psi_2^{22}+\psi_1^{21}+\psi_1^{12}$}\label{01case5}
We will begin by computing the coboundary operator with representatives from the $\Cn{n}{2}$ and $\Cn{n}{1}$ spaces:
\begin{align*}
 D(\ph^I_2)=&\ph^{I2}_2+(-1)^{I+1}\ph^{2I}_2+\ph^{I1}_1+(-1)^{I+1}\ph^{1I}_1+(-1)^I\ph^I_2\psi^{22}_2 \\ &+(-1)^I\ph^I_2\psi^{21}_1+(-1)^I\ph^I_2\psi^{12}_1 \\
 D(\ph^I_1)=&(-1)^{I+1}\ph^{2I}_1+\ph^{I2}_1+(-1)^I\ph^I_1\psi^{22}_2+(-1)^I\ph^I_1\psi^{21}_1 \\ &+(-1)^I\ph^I_1\psi^{12}_1
\end{align*}
We decompose our spaces in the following manner,
\begin{align*}
 \Cn{n}{a,k}=&\langle\ph^{2^k1I}_a\mid\ell(I)=n-k-1\rangle &\Pn{n}a=&\langle\ph^{2^n}_a\rangle \\
 \Cn{n}{a,0,k}=&\langle\ph^{1^k2I}_a\mid\ell(I)=n-k-1\rangle &\Sn{n}a=&\langle\ph^{1^n}_a\rangle \\
\end{align*}
With these defintions there is ambiguity when $n=0$, so we define $\Cn{0}a=\Sn{0}a$. Note that we can write $\Cn{n}{a,0}$ as
\begin{equation*}
 \Cn{n}{a,0}=\Sn{n}a \oplus \bigoplus_{k=1}^{n-1} \Cn{n}{a,0,k}
\end{equation*}

 We now decompose the $D$ operator as follows,
\begin{align*}
 D=D'_2+D''_2+D^k_1+D^0_1:\Cn{n}{2,k}\ra&\Cn{n+1}{2,k}\oplus\Cn{n+1}{2,k+1}\oplus\Cn{n+1}{1,k}\oplus\Cn{n+1}{1,0} \\
 D=D'_1+D''_1:\Cn{n}{1,k}\ra&\Cn{n+1}{1,k}\oplus\Cn{n+1}{1,k+1}
\end{align*}
For convenience we denote $D_1=D'_1+D''_1$ and $D_2=D'_2+D''_2$. With this decomposition we obtain the following relations,
\begin{align*}
 (D'_2)^2=&0 & (D'_1)^2=& 0 &   D_1D^k_1=&-D^k_1D_2 \\
 D'_2D''_2=&-D''_2D'_2 & D'_1D''_1=&-D''_1D'_1  & D_1D^0_1=&-D^0_1D_2 \\
 (D''_2)^2=&0 & (D''_1)^2=&0 & &
\end{align*}

Some sort of transitional sentence that will start the proving process.
\begin{align*}
 D''_a(\ph^{2^k1I}_a)=&\left\{\begin{array}{ll}
                               0, & k\text{ even} \\
                               \pm\ph^{2^{k+1}1I}_a, & k \text{ odd}
                              \end{array}
                              \right.
&D(\ph^{2^n}_a)=&\left\{\begin{array}{ll}
                         0, & n\text{ even} \\
                         \pm\ph^{2^{n+1}}_a, & n \text{ odd}
                        \end{array}
                        \right.
\end{align*}

Let $\ph\in\Cn{n}2$ and $\xi\in\Cn{n}1$ satisfy $\ph+\xi\in\ker{D}$. We can write $\ph$ and $\xi$ as
\begin{align*}
 \ph=&\ph_0+\dots+\ph_{n-1}+\ph_P \\
 \xi=&\xi_0+\dots+\xi_{n-1}+\xi_P.
\end{align*}
Where $\xi_k\in\Cn{n}{1,k}$ and $\xi_P\in\Pn{n}1$. Since $D(\ph+\xi)=0$ we see,
\begin{align}
D_2(\ph)=&0 \nonumber \\
D_1^k(\ph_0)+D_1^0(\ph)+D_1(\xi_0)=&0,\quad \text{if }k=0  \label{eqn-xi0}\\
\quad D_1^k(\ph_k)+D'_1(\xi_k)+D''_1(\xi_{k-1})=&0,\quad \text{If }k\ge1.\label{eqn-xik}
\end{align}
Our goal is to find $\alpha\in\Cn{n-1}2$ and $\eta\in\Cn{n-1}1$ such that $D(\alpha+\eta)=\ph+\xi$. Note that we can write $\alpha$ and $\eta$ as a sum of terms as we did with $\ph$ and $\xi$. Using THEOREM on $D_2$ we can find $\alpha_{k\ge1}\in\Cn{n-1}{2,k\ge1}$ such that $D_2(\alpha_{k\ge1})=\ph_{k\ge1}$. Now we find $\eta_{k\ge1}$ to satisfy (\ref{eqn-xik}). For $k=n+1$ we have
\begin{equation*}
 D''_1(\xi_n)=0
\end{equation*}
and since $H(D''_1)=0$ on the subcomplex $\Cn{n}{1,k\ge1}$, we can find an $\eta_{n-1}\in\Cn{n-1}1$ such that $D''_1(\eta_{n-1})=\xi_n$. Assume we have shown for $k+1$, or assume we have shown
\begin{equation*}
 \xi_{k+1}=D^k_1(\alpha_{k+1})+D'_1(\eta_{k+1})+D''_1(\eta_k).
\end{equation*}
Then we have
\begin{align*}
 0=&D^k_1(\ph_{k+1})+D'_1(\xi_{k+1})+D''_1(\xi_k) \\
  =&D_1^kD'_2(\alpha_{k+1})+D^k_1D''_2(\alpha_k)+D'_1D^k_1(\alpha_{k+1}) \\
  &+D'_1D''_1(\eta_k)+D''_1(\xi_k) \\
  =&D''_1(\xi-D_1^k(\alpha_k)-D'_1(\eta_k))
\end{align*}
Thus $\xi_k=D^k_1(\alpha_k)+D'_1(\eta_k)+D''_1(\eta_{k-1})$ for some $\eta_{k-1}\in\Cn{n-1}{1,k-1}$. This holds until $k=1$ in which case we show
\begin{align*}
 0=&D^k_1(\ph_2)+D'_1(\xi_2)+D''_1(\xi_1) \\
  =&D^k_1D''_2(\alpha_1)+D'_1D''_1(\eta_1)+D''_1(\xi_1) \\
  =&D''_1(\xi_1-D^k_1(\alpha_1)-D'_1(\eta_1)).
\end{align*}
However, $D''_1$ is injective when $k=1$ so $\xi_1=D^k_1(\alpha_1)+D'_1(\eta_1)$.

When $k=0$ a separate technique is required. We begin by showing that $\la^{1}$ commutes with the $D_a$ operators.
\begin{align*}
 D_2(\la^1\ph^I_2)=&\ph^{1I2}_2+(-1)^{I}\ph^{21I}_2+(-1)^I\ph^I_2\psi^{22}_2+(-1)^{I+1}\ph^{21I}_2 \\
 &+(-1)^I\la^1\ph^I_2\psi^{21}_1+(-1)^{I+1}\ph^{12I}_2+(-1)^I\la^1\ph^I_2\psi^{12}_1  \\
 =&\la^1D_2(\ph^I_2)
\end{align*}
The same proof will hold for $D_1\la^1=\la^1D_1$. This implies that $\la^{1^k}$ commutes with $D_a$ and notice that
\begin{equation*}
 \la^{1^k}:\Cn{n}a\mapsiso\Cn{n}{a,0,k}
\end{equation*}

Similar to before we can write $\ph_0$ and $\xi_0$ as,
\begin{align*}
 \ph_0=&\ph_{0,2}+\dots+\ph_{0,n-1}+\ph_{0,S} \\
 \xi_0=&\xi_{0,2}+\dots+\xi_{0,n-1}+\xi_{0,S}
\end{align*}

Let $\ph'_k\in\Cn{n-k}{2,k}$ be such that $\la^{1^k}\ph'_k=\ph_{0,k}$. Now we rewrite our condition on $\ph_0$ as,
\begin{equation*}
 D_2(\ph_0)=0 \Rightarrow 0=D_2(\ph_{0,k})=\max_k\la^{1^k}D_2(\ph'_k).
\end{equation*}
Since we pulled out the maximum number of $1$'s, then if $\ph'_k\not\in\Cn02$ then $\ph'_k$ has a leading $2$, and since $H(D_2)=0$ on $\Cn{n}{2,k\ge1}$ and $\Pn{n}2$ we can conlclude that $\ph'_k=D_2(\alpha'_k)$ for some $\alpha'_k\in\Cn{n-k}{2,k}$ and thus,
\begin{equation*}
 \ph_0=D(\alpha_0).
\end{equation*}
This means we are left only with the term $\ph_{0,S}$ from $\ph$.

Now consider $\xi_0$, define $\xi'_k\in\Cn{n-k}{1,k}$ by $\xi_{0,k}=\la^{1^k}\xi'_k$. We want to find $\eta_{0,k}\in\Cn{n-1}{1,0,k}$ such that
\begin{align*}
\xi_{0,1}=&D^0_1(\alpha_{0,k\ge1})+D_1(\eta_{0,1})+D^k_1(\alpha_{0,1}) \\
\xi_{0,k\ge2}=&D^0_1(\alpha_{0,k-1})+D_1(\eta_{0,k})+D^k_1(\alpha_{0,k})
\end{align*}
However, we've already chosen our $\alpha$'s, so we really want
\begin{align*}
 D_1(\eta_{0,1})=&\xi_{0,1}-D^0_1(\alpha_{0,k\ge1})-D^k_1(\alpha_{0,1}) \\
 D_1(\eta_{0,k})=&\xi_{0,k\ge2}-D^0_1(\alpha_{0,k-1})-D^k_1(\alpha_{0,k})
\end{align*}

For $k\ge2$ we have
\begin{align*}
 0=&D^0_1(\ph_{0,k-1})+D^k_1(\ph_{0.k})+D_1(\xi_{0,k}) \\
  =&D_1^0D_2(\alpha_{0,k-1})+D^k_1D_2(\alpha_{0,k})+D_1(\xi_{0,k}) \\
  =&D_1(\xi_{0,k}-D^0_1(\alpha_{0,k-1})-D^k_1(\alpha_{0,k})) \\
  =&\max_k\la^{1^k}D_1(\xi'_k-D^0_1(\alpha'_{k-1})-D^k_1(\alpha'_k)).
\end{align*}
By an argument similar to the $\ph_0$, if $\xi'_0\not\in\Cn{0}1$ then we can find an $\eta'_k\in\Cn{n-k}{1,k}$ such that $\xi'_k=D^0_1(\alpha'_{k-1})+D^k_1(\alpha'_k)+D_1(\eta'_k)$, more specifically if $k\ge2$ we can find $\eta_{0,k\ge2}$ such that
\begin{equation*}
 \xi_{0,k\ge2}=D^0_1(\alpha_{0,k-1})+D^k_1(\alpha_{0,k})+D_1(\eta_{0,k})
\end{equation*}

If $k=1$ then we have
\begin{align*}
 0=&D^k_1(\ph_{0,1})+D^0_1(\ph_{k\ge2})+D_1(\xi_{0,1}) \\
  =&D^k_1D_2(\alpha_{0,1})+D_1^0D_2(\alpha_{k\ge2})+D_1(\xi_{0,1}) \\
  =&-D_1D^k_1(\alpha_{0,1})-D_1D^0_1(\alpha_{k\ge2})+D_1(\xi_{0,1})
  =&D_1(\xi_{0,1}-D^k_1(\alpha_{0,1})-D^0_1(\alpha_{k\ge2}))
\end{align*}
Thus we can write
\begin{equation*}
 \xi_{0,1}=D^0_1(\alpha_{0,k\ge2})+D_1(\eta_{0,1})+D^k_1(\alpha_{0,1})
\end{equation*}
for some $\eta_{0,1}\in\Cn{n-1}{2,0,1}.$

Finally consider the spaces $\Sn{n}2$ and $\Sn{n}1$. To begin we need to be more specific about the action of the coboundary operator on these spaces.
\begin{align*}
 D(\ph^{1^n}_2)=&\left\{\begin{array}{ll}
                        2\ph^{1^{n+1}}_1, & n\text{ odd} \\
                        0, & n\text{ even}
                       \end{array}\right. \\
 D(\ph^{1^n}_1)=&0
\end{align*}
When $n$ is odd $D:\Sn{n}2\mapsiso \Sn{n+1}1$. Thus we obtain the diagram in Figure \ref{fc5}. The spaces which are boxed are isolated spaces that go to zero, or contribute cohomology. Thus we see the cohomology is given by
$$H^n=\left\{\begin{array}{ll}
              \psi_2,\psi_1, & n=0 \\
              \phi^{1^n}_1, & n \text{ odd} \\
              \phi^{1^n}_2, & n\ne0 \text{ and }n\text{ even}
             \end{array}\right.$$
\begin{figure}[h!]
 $$\xymatrix{
 \Sn{0}2 \ar[rd]^0 & \Sn{0}1 \ar[d]^0
 \save [l].[]*[F]\frm{}\restore\\
 \Sn{2}2 \ar[rd]^\thicksim & *+[F]{\Sn{2}1} \ar[d]^0 \\
 *+[F]{\Sn{1}2} \ar[rd]^0 & \Sn{1}1 \ar[d]^0 \\
 \Sn{3}2 \ar[rd]^\thicksim & *+[F]{\Sn{3}2} \ar[d]^0 \\
 & \vdots
 }$$
 \caption{}\label{fc5}
\end{figure}


\subsection{$d_6=\psi_1^{22}$}\label{01case6}
For case six, the methods employed in the previous cases completely fail, and
as a result, a different method is necessary. Once again, we begin by computing
the bracket of \emph{$d$} with a general element in $\Cn{n}{2}$ and $\Cn{n}{1}$.
\begin{align*}
D(\ph^I_2)&=\ph^{I2}_1+(-1)^{I+1}\ph^{2I}_1+(-1)^I\ph^I_2\ph^{22}_1\\
D(\ph^I_1)&=(-1)^{I}\ph^I_1\ph^{22}_1
\end{align*}
Therefore, we have decompositions
\begin{align*}
D=D_2+D_1&:\Cn{n}{2}\ra \Cn{n+1}{2}\oplus \Cn{n+1}{1}\\
D&:\Cn{n}{1}\ra \Cn{n+1}{1}
\end{align*}
Note that $D_1$ is injective and that
\begin{equation*}
D_2^2=0,\qquad D_1D_2=-DD_1.
\end{equation*}

Thus $D_2$ is a coboundary operator on $\Cn{n}{2}$, giving a cohomology
$H^n_2=H^n(D_2)$. We first compute this cohomology, and use it to compute the
cohomology in general. Define the Decleene map
$\theta=\lambda^{21}+\lambda^{12}$. Then we claim that $\theta$ commutes with
$D_2$ on $\Cn{n}{2}$, and with $D$ on $\Cn{n}{1}$. To see this, note that
\begin{align*}
D_2\theta(\pha I2)&=D_2(\pha{21I}2+\pha{12I}2)\\
=&\s{I+1}\ph^{222I}_2+\s{I}\ph^{222I}_2\s{I}\lambda^{21}\ph^I_2\psi^{22}_1+\s{I}\pha{I}2\psa{22}2 \\
=&\theta D_2(\pha{I}2)
\end{align*}
The proof that the DeCleene map commutes with $D$ on $\Cn{n}{1}$ is similar. In fact,
note that the action of $D$ on $\Cn{n}{1}$ is essentially the same as $D_2$ on $\Cn{n}{2}$.

Next, note that if $\ph$ is a $D_2$-coboundary,  then every term in $\ph$ must
have a double $2$. Therefore, any $D_2$-cocycle which has a term without a
double $2$ must be nontrivial.  In particular, the 0-cocycle $\ph_2$ and the
2-cocycle $\pha 22$ are nontrivial $D_2$-cocycles.  Define the Decleene
cocycle $\dec^n_2$ and $\dec^n_1$ by
\begin{align*}
\dec^{2n}_2&=\theta^n\ph_2, &\dec^{2n+1}_2&=\theta^n\pha 22,\\
\dec^{2n}_1&=\theta^n\ph_1, &\dec^{2n+1}_1&=\theta^n\pha 21.
\end{align*}
Then $\dec^n_2$ is a nontrivial $D_2$-cocycle. Also $\dec^n_1$ is nontrivial if
we consider only the cohomology of $D$ restricted only to the $F$ space. We
shall discuss later when it is a nontrivial cocycle on the whole space
$C^n=\Cn{n}{1}\oplus \Cn{n}{2}$.

Let $B^n_2$ be the space of $D_2$ $n$-coboundaries, $Z^n_2$ be the
$n$-cocycles, $z_n=\dim(Z^n_2)$, $b_n=\dim(B^n_2)$ and $h_n=\dim H^n_2$. Then
$h_n=z_n-b_n$ and $z_n+b_{n+1}=2^n$. Because there is a nontrivial Decleene
cocycle in every degree, we know that $h_n\ge 1$. We wish to show that $h_n=1$.

To see this, note that $\theta_n:B^n\ra B^{n+1}$ is injective,
$D_2\circ\lambda^1:\Cn{n}{2}\ra B^{n+1}$ is also injective, and the images of these
operators are independent subspaces.  As a consequence, we must have
$b_{n+1}\ge b_n+2^n$. We will show that $b_n+b_{n+1}\ge 2^n-1$ for $n\ge 0$.
Since $b_0=b_1=0$, and $b_1=1$ by direct computation, the formula holds for
$n\le 2$. Suppose it holds for $k<n$, and that $n\ge 2$. Then
\begin{equation*}
b_{n}+b_{n+1}\ge b_{n-2}+2^{n-2}+b_{n-1}+2^{n-1}\ge 2^{n-2}+2^{n-2}+2^{n-1}=2^{n}.
\end{equation*}
Thus, by induction, the formula  holds for all $n$. Using this formula, we
obtain
\begin{align*}
1\le h_n=z_n-b_n=2^n-b_{n+1}-b_n\le 2^n-(2^n-1)=2.
\end{align*}

First, let us say that $\dec^n_2$ extends to a $D$-cocycle if there is some
$\eta\in \Cn{n}{1}$ such  that $\dec^n_2+\eta$ is a $D$-cocycle. If $\dec^n_2$
extends, then let $\dec^n=\dec^n_2+\eta$ be some arbitrary extension of
$\dec^n_2$.

Suppose that $D(\ph+\xi)=0$ for some $\ph\in \Cn{n}{2}$ and $\xi\in \Cn{n}{1}$. Then
$D_1(\ph)+D(\xi)=0$ and $D_2(\ph)=0$. In fact, the second equation follows from
the first one. For, suppose the first equality holds. Then
\begin{equation*}
D_1D_2(\ph)=-DD_1(\ph)=D^1(\xi)=0.
\end{equation*}
Using the fact that $D_1$ is injective, we see that $D_2(\ph)=0$. Now we can
write $\ph=a\dec^n_2+D_2(\alpha)$ for some $\alpha\in \Cn{n-1}{2}$, because we know
that $h_n=1$.

Note that if $\dec^n_2$ does not extend to a $D$-cocycle, then $a=0$. This is
because
\begin{align*}
0&=D_1(\ph)+D(\xi)\\&=D_1(a\dec^n_2)+D_1D_2(\alpha)+D(\xi)\\&=D_1(a\dec^n_2)-DD_1(\alpha)+D(\xi)
\\&=D_1(a\dec^n_2)-D(D_1(\alpha)-\xi)),
\end{align*}
so that if $a\ne0$ we have $D_1(\dec^n_2)=D(\eta)$, where
$\eta=\tfrac1a(D_1(\alpha)-\xi)$. Next, we claim  that
$\ph+\xi=b\dec^n_1+D(\alpha+\beta)$ for some $\beta\in \Cn{n-1}{1}$. To see this,
first suppose that $a=0$. Then
\begin{align*}
0=D(\ph+\xi)=D_1D_2(\alpha)+D(\xi)=-D(D_1(\alpha)+D(\xi)=D(\xi-D_1(\alpha)).
\end{align*}
Thus $\xi-D_1(\alpha)$ is a $D$-cocycle lying in $\Cn{n}{1}$, which means it can be
written in the form $\xi-D_1(\alpha)=b\dec^n_1+D(\beta)$ for some $\beta\in
\Cn{n-1}{1}$. But this means $\ph+\xi=b\dec^n_1+D(\alpha+\beta)$, as desired.

On the other hand, if $a\ne0$ then $a\dec^n_2=a\dec^n+a\eta$, where $\eta\in
\Cn{n-1}{1}$, so $\ph=a\dec^n+a\eta+D_2(\alpha)$, and then
\begin{equation*}
0=D(\ph+\xi)=D(a\eta)+D_1D_2(\alpha)+D(\xi)=D(\xi+a\eta-D_1(\alpha)).
\end{equation*}
Thus in this case, we can express $\xi+a\eta-D_1(\alpha)=b\dec^n_1+D(\beta)$,
so we obtain
\begin{equation*}
\ph+\xi=a\dec^n+a\eta+\xi=a\dec^n+b\dec^n_1+D(\alpha+\beta).
\end{equation*}
From the equation above, it follows that the dimension of $H^n$ is at most 1,
depending on whether $\dec^n_2$ extends to a $D$-cocycle and whether $\dec^n_1$
is a nontrivial cocycle.

Now we show that the non triviality of $\Ch^n_f$ is linked to the whether or
not we can extend $\Ch^{n-1}_2$.

Suppose $\dec^n_1$ is trivial, ie. $\dec^n_1=D(\ph+\xi)$ for some $\ph\in
\Cn{n-1}{2}$, $\xi\in \Cn{n-1}{1}$. $D_2(\ph)=0$ so $\ph=a\dec^{n-1}_2+D_2(\alpha)$ for
some $\alpha\in \Cn{n-1}{2}$. If $\Ch^{n-1}_2$ extends, so
$\dec^{n-1}_2=\dec^{n-1}+\eta$ then
\begin{align*}
\dec^n_1&=D_1(\ph)+D(\xi)=aD(\eta)+D_1D_2(\alpha)+D(\xi)\\&=
D(\eta)-DD_1(\alpha)+D(\xi)=D(\xi+\eta-D_1(\alpha)).
\end{align*}
But then $\dec^n_1$ is a coboundary in the $\Cn{n}1$ space, which is impossible. Thus
if $\dec^n_1$ is trivial, $\Ch^{n-1}_2$ does not extend to a $D$-cocycle.

On  the other hand, suppose $\dec^{n-1}_2$ does not extend to a $D$-cocycle.
Then $D_1(\dec^{n-1}_2)$ is a $D$-cocycle, lying in $\Cn n1$, which is nontrivial in
terms of the $D$-cohomology restricted to $\Cn{n}1$, so we must have
$D_1(\dec^{n-1}_2)=a\dec_1^n+D(\beta)$ for some $\beta\in \Cn{n-1}{1}$, where
$a\ne0$. But then $\dec^n_1=D(\tfrac1a(\dec^{n-1}_2-\beta))$, and therefore
$\dec^n_1$ is trivial.

We now show that the $\dec_1^n$ is non-trivial when $n$ is both even and odd. First when $n=2\mod1$ suppose that $\dec_1^n=D(\ph+\xi)$, where $\ph\in\Cn{n-1}{2}$ and $\xi\in\Cn{n-1}{1}$. Since we also want $\ph$ to be a $D_2$ cocylce, we conclude that $\ph=a\dec_2^{n-1}+D_2(\alpha)$ for some $\alpha\in\Cn{n-1}{2}.$ Since $\dec^n_1$ contains terms of the form $\ph^{(21)^{\frac{n-1}2}2}_1$ then, by assumption, we must have this term appear in $D_1(\ph)=aD_1(\dec^{n-1}_2)+D_1D_2(\alpha)$. Clearly this term won't appear in $D_1D_2(\alpha)$, and $D_1(\dec^{n-1}_2)$ will produce two terms of the proper form but of opposite sign. So $\dec^{n}_1$ must be non-trivial.

Now assume $n=0\mod2$. Following the same lines as the previous arguement, we see $\dec^n_1$ will contain terms of the form $\ph^{(21)^\frac n2}_1$ and $\ph^{(12)^\frac n2}_1$. But $\dec^{n-1}_2$, will have no terms with a $1$ on the edge and no double $2$. Thus $\dec^n_1$ must be non-trivial.

We are able to conclude that $h^n(d)=2$ for all $n$.

\section{Infinitesimal Deformations}
To compute the infinitesimal deformations we only need consider cohomology in degree two. That leaves only $d_2$, $d_5$ and $d_6$; each of which will be considered separately.
\subsection{$d_2$}
The cohomology of $d_2$ is given by $\psa{11}1$. We determine that
$$d_t=\psa{22}2+t\psa{11}1$$
is isomorphic to $d_1$ when $t\ne0$.
\subsection{$d_5$}
The cohomology of $d_5$ is given by $\psa{11}2$. We determine that
$$d_t=\psa{22}2+\psa{21}1+\psa{12}1+t\psa{11}2$$
is isomorphic to $d_1$ when $t\ne0$.
\subsection{$d_6$}
The cohomology of $d_6$ is given by $\psa{21}2+\psa{12}2+\psa{11}1$ and $\psa{21}1+\psa{12}1$. We determine that
$$d_t=\psa{22}1+t_2(\psa{21}2+\psa{12}2+\psa{11}1)+t_1(\psa{21}1+\psa{12}1)$$
is isomorphic to $d_2$ when $t_1^1=t_2$, to $d_2$ when $t_2^1=\tfrac{3t_1^1}4$, and when $t_2$ and $t_1$ are not both zero to $d_1$.
\section{Versal Deformations}
In this case the versal deformations coincide exactly with the infinitesimal deformations. This can be seen be taking the bracket of the infinitesimal deformation with itself. If the result is zero then the deformations coincide.
\subsection{Diagram of Deformations}
\begin{figure}[h!]
 $$\xymatrix@R=5mm@C=3mm{
  &&& \bullet^{d_6}\ar@{~>}[rd]\ar@{~>}[ld]\\
  \bullet^{d_3} & \bullet^{d_4} & ^{d_2}\bullet \ar@{~>}[rd] & & \bullet^{d_5} \ar@{~>}[ld] \\
  &&& \bullet_{d_1}&
  }$$
  \caption{Infinitesimal and Versal deformations of an $0|2$-dimensional vector space.}\label{dia}
\end{figure}
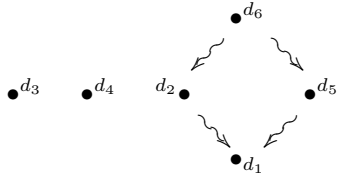

For a visual of the deformations we provide Figure \ref{dia}.


\bibliographystyle{amsplain}
\providecommand{\bysame}{\leavevmode\hbox to3em{\hrulefill}\thinspace}
\providecommand{\MR}{\relax\ifhmode\unskip\space\fi MR }
\providecommand{\MRhref}[2]{%
  \href{http://www.ams.org/mathscinet-getitem?mr=#1}{#2}
}
\providecommand{\href}[2]{#2}

\end{document}